\documentclass{amsart}%
\usepackage{amsmath}
\usepackage{amsfonts}
\usepackage{amsthm}
\usepackage{amssymb}
\newtheorem{theorem}{Theorem}
\newtheorem*{theorem*}{Theorem}
\newtheorem*{acknowledgement*}{Acknowledgement}

\newtheorem{example}[theorem]{Example}

\newtheorem{lemma}[theorem]{Lemma}

\newtheorem{proposition}[theorem]{Proposition}
\newtheorem{remark}[theorem]{Remark}

\newcommand{\Mnp}[0]{\mathcal{M}^{n+1}}

\newcommand{\R}[1]{{\mathbb{R}}^{#1}}

\newcommand{\hess}[1]{\nabla\nabla #1}

\newcommand{\gt}[0]{\tilde{g}}
\newcommand{\om}[0]{\omega}
\newcommand{\op}[0]{\omega^{\prime}}
\newcommand{\opp}[0]{\omega^{\prime\prime}}
\newcommand{\fp}[0]{f^{\prime}}
\newcommand{\fpp}[0]{f^{\prime\prime}}
\newcommand{\dt}[0]{\frac{d}{dt}}
\newcommand{\ddt}[1]{\frac{d #1}{dt}}
\newcommand{\ddtt}[1]{\frac{d^{2}#1}{dt^2}}
\newcommand{\ddttt}[1]{\frac{d^{3} #1}{dt^3}}
\newcommand{\delt}[0]{\widetilde{\nabla}}

\begin{document}
\title[On rotationally invariant shrinking solitons]{On rotationally invariant shrinking gradient Ricci solitons}

\author[B. Kotschwar] {Brett Kotschwar}
\address{University of California, San Diego}

\email{bkotschw@math.ucsd.edu}
\begin{abstract}
  In this paper we study the gradient Ricci shrinking soliton equation on 
  rotationally symmetric manifolds of dimension 
  three and higher and prove that the only complete examples of such metrics
  on $S^n$, $\R{n}$ and $\R{}\times S^{n-1}$ are, respectively, the round, flat, 
  and standard cylindrical metrics.
\end{abstract}
\maketitle
\section{Introduction}
Recall that a Riemannian manifold $(\Mnp, g)$ 
is said to be a \emph{shrinking gradient Ricci soliton} provided
\begin{equation}
\label{eq:sol}
	Rc(g) + \hess f - \lambda g = 0
\end{equation}
for some $f\in C^{\infty}(\Mnp)$ and $\lambda > 0$. 
The analogous objects in the cases $\lambda = 0$ and $\lambda <0$ are known as
\emph{steady} and \emph{expanding} solitons, respectively, and 
the triple $(g, f, \lambda)$ is referred to 
as a \emph{(gradient) soliton structure} on the manifold $\Mnp$.
Solitons correspond to self-similar solutions of the Ricci flow and arise 
commonly in the analysis of its singularities.

After the class of Einstein metrics, perhaps the 
most natural place to look for solitons is among the 
rotationally invariant metrics, and, in all dimensions greater than one, this class has been shown
to contain complete, non-trivial examples of steady and expanding solitons 
(cf., e.g, \cite{I} and \cite{RFV2}).
However, one does not expect to find corresponding examples in the shrinking 
case. The purpose of the present work is to confirm this expectation. 
Our main result is the following
\begin{theorem}
\label{thm:coverclass}
Suppose $n \geq 2$, and $(g, f, \lambda)$ is a complete, rotationally invariant shrinking soliton 
structure on a on a manifold $\Mnp$ diffeomorphic to one of $S^{n+1}$, $\R{n+1}$ or
 $\R{}\times S^{n}$.  Then,
\begin{enumerate}
	\item \label{it:sn}
	 If $\Mnp\approx S^{n+1}$, $g$ is isometric to a round sphere
		and $f\equiv \mathrm{const}$.
	\item \label{it:rn}
		If $\Mnp\approx \R{n+1}$, $g$ is flat.
	\item \label{it:cyl}
		If $\Mnp\approx \R{}\times S^n$, 
		$g$ is isometric
		to the standard cylinder $dr^2 + \omega_0^2g_{S^{n}}$
		 of radius 
		$\omega_0 =\sqrt{(n-1)/\lambda}$
		and $f= f(r) = (n-1)r^2/(2\omega_0^2) + \mathrm{linear}$. 
\end{enumerate}
\end{theorem}
This theorem is known (in greater generality) in dimensions less than four.
In dimension two, work of Hamilton (\cite{H}, \cite{H2}) shows that the only 
complete shrinking solitons are the flat metric on $\R{2}$ and the
round metric on $S^2$.  In fact, as Hamilton observes in \cite{H},
surface solitons are necessarily rotationally symmetric--an observation
used, for example, by Chen, Lu, and Tian 
in their note \cite{CLT} to provide an alternative proof
of the uniqueness of the constant
curvature soliton on $S^2$.
In a sense, then, our result is an extension of these 
two-dimensional findings to higher dimensions.

In dimension three, gradient solitons need not be rotationally symmetric,
however Perelman \cite{P2} has shown that the only 
complete examples of shrinking solitons with non-negative sectional 
curvature are the flat metric on $\R{3}$, the round metric on $S^3$, and the
standard metric on the cylinder $\R{}\times S^{2}$.  
As the Hamilton-Ivey estimate (\cite{H}, \cite{I}) 
implies that three-dimensional ancient solutions
are necessarily of non-negative curvature, 
Perelman's argument classifies all complete three-dimensional 
shrinking solitons
  
In higher dimensions, a consequence of the recent work of B\"{o}hm and 
Wilking \cite{BW}, is that the only
compact shrinking soliton with $2$-positive curvature operator 
is the round sphere.
In the K\"{a}hler category, Feldman, Ilmanen, and Knopf
\cite{FIK} have obtained results similar to our own 
for K\"{a}hler-Ricci shrinking solitons
under the assumption of $U(n)$-invariance. 
In particular, their Proposition 9.2--
that the flat metric is the only complete 
$U(n)$-invariant gradient shrinking soliton on ${\mathbb{C}}^n$--is 
the K\"{a}hler analog 
of the second case of our Theorem \ref{thm:coverclass}. Their paper
also provides non-trivial examples of $U(n)$-invariant 
gradient solitons (of all types) on other spaces.  

To prove Theorem \ref{thm:coverclass}, 
we show first that, under the constraint (\ref{eq:sol}),
 the rotational symmetry of a non-flat metric implies the
rotational symmetry of the gradient function.  
This reduces the proof of the proposition to the study
of a certain second order system of non-linear ODEs.  
By a change of variables
due to Robert Bryant and Tom Ivey \cite{RFV2}, we are able to further reduce the 
problem to the study of an equivalent first-order system 
amenable to phase-plane analysis.   
For the case $\Mnp\approx S^{n+1}$, we show that any candidate
metric must have positive curvature operator and then apply
the result of B\"{o}hm and Wilking. It is no doubt possible to prove this
case solely from the analysis of the system of ODE (without using
any of the dynamic properties of the Ricci flow), however, since we feel the
result may not come as a surprise to experts,
we content ourselves here with a proof by the most readily available means. 
For the non-compact cases, we use the 
criterion of completeness to eliminate all but the two 
standard metrics by their asymptotic behavior.

\section{The soliton condition for rotationally \\ symmetric manifolds.}

\subsection{The warped-product metric}
Fix  $n > 1$ and let $\gt$ denote the metric on $S^n$ of 
constant sectional curvature $1$.
For
$-\infty \leq A < \Omega \leq \infty$, and positive 
functions $\omega\in C^{\infty}(A,\Omega)$, we shall consider the warped-product metric 
$g= dr^2 + \omega^2(r)\gt$ on 
the cylinder ${\mathcal{C}}_{A,\Omega}:= (A, \Omega)\times S^{n}$. 
For this metric, we have the following standard 
\begin{lemma}
\label{lem:smooth}
  The metric $g= dr^2 + \om(r)^2\gt$ on ${\mathcal{C}_{0, \Omega}}$ extends to a smooth metric on $B_{\Omega}(\mathbf{0})\subset R^{n+1}$
if and only if $\Omega > 0$ and
\[
	\lim_{r\to 0}\om(r) = 0\mbox{,} \,\,\lim_{r\to 0}\op = 1 \mbox{, and}\,\,
\lim_{r\to 0} \frac{d^{2k}\om}{dr^{2k}}(r) = 0
\,\,\mbox{for all}\, k.
 \]
 
 The metric extends to a smooth metric on $S^{n+1}$ if and only if 
in addition $\Omega < \infty$ and
\[
	\lim_{r\to \Omega}\om(r) = 0\mbox{,} \,\,\lim_{r\to \Omega}\op = - 1 \mbox{, and}\,\,
\lim_{r\to \Omega} \frac{d^{2k}\om}{dr^{2k}}(r) = 0
\,\,\mbox{for all}\, k.
 \]
\end{lemma}
\begin{proof}
  See, e.g., \cite{CK}, Lemma 2.10.
\end{proof}

Thus we may accomplish the proof of Theorem \ref{thm:coverclass}
by the study of soliton structures on $\mathcal{C}_{A,\Omega}$.  We begin by recording the expressions
of a few geometric quantities associated to $g$.

The sectional curvatures of $g$
are
\begin{equation}
\label{eq:sect}
	\nu_1 = -\frac{\opp}{\om}\,\, \mbox{  and  }\,\, \nu_2 = \frac{1-(\op)^2}{\om^2},
\end{equation}
where $\nu_1$ denotes the curvature of two-planes tangent to the radial direction,  
and $\nu_2$ the curvature of orbital two-planes.
The Ricci curvature of $g$ has the form
\begin{equation}
\label{eq:rc}
	Rc(g) = -n\frac{\opp}{\om}dr^2 
		+ \left((n-1)(1-(\op)^2) -\om\opp\right)\gt
\end{equation}
where the prime denotes differentiation with respect to $r$. 

A routine computation also shows that 
for a smooth function $f$ on ${\mathcal{C}}_{A,\Omega}$, one has
\begin{equation*}
\label{eq:fhess}
	\hess{f} = \begin{cases}
			\nabla_0\nabla_0 f &= f_{00}\\
			\nabla_0\nabla_i f &= f_{i0} - (\omega^{\prime}/\om)
			f_i\\
			\nabla_i\nabla_j f &= 
		\widetilde{\nabla}_i\widetilde{\nabla}_j f + \om\om^{\prime}f_0\gt_{ij}	
		\end{cases}
\end{equation*}
in local coordinates $(\theta^0 = r,\theta^i)$.\footnote{ Here and throughout, when working
in coordinates $\theta^0,\ldots,\theta^n$, we shall use roman letters to denote
indices $1,\ldots,n$ and use a tilde to denote quantities (Levi-Civita connection, 
curvature, etc.) associated
to the metric $\gt$ on $S^n$.  In particular, 
$\widetilde{\nabla}_i\widetilde{\nabla}_j f$  represents the hessian of $f(r,\cdot)$ considered as a function on $S^{n}$.}  

In these coordinates, (\ref{eq:sol}) has the expression
\begin{align}
\label{eq:rsol00}
    &f_{00} - n\frac{\opp}{\om}-\lambda =0\\
\label{eq:rsol0i}
  &f_{0i} -\frac{\op}{\om}f_i =0\\
\label{eq:rsolij}
  &\widetilde{\nabla}_i\widetilde{\nabla}_j f 
  +\left[(n-1)(1-(\op)^2)-\om\opp- \om\om^{\prime}f_0-\lambda\om^2\right]\gt_{ij} =0.
\end{align}	

The above system involves the partial derivatives of $f$, and despite the rotational symmetry 
of the metric and the Ricci curvature, there is no \textit{a priori}
reason to assume that $f$ shares this symmetry.  
However, as we show next, this symmetry is implied unless $g$ is flat,
and thus for the proof of Theorem \ref{thm:coverclass}, there is no loss in 
restricting our attention to the case $f= f(r)$.  

When $f$ is a radial function, equations
(\ref{eq:rsol00})-(\ref{eq:rsolij}) reduce to
\begin{align}
\label{eq:rsys}
	\begin{cases}
	\fpp -\lambda &= n\opp/\om\\
	\om\op\fp -\lambda\om^2 &= \om\om^{\prime\prime} + (n-1)
\left((\op)^2 -1\right).
	\end{cases}
\end{align}

\begin{proposition}
  Suppose that $(g,f,\lambda)$ is a complete, rotationally symmetric gradient shrinking
soliton structure
  on $\Mnp\approx S^{n+1}$, $\R{n+1}$, or $\R{}\times S^n$.

If the function $f$ is not rotationally symmetric, then  
$\Mnp\cong \R{n+1}$ and $g$ is the flat metric.
\end{proposition}
\begin{proof}
  Write $g = dr^2 + \om^2(r)\gt$ for $r\in (A,\Omega)$, 
  and fix local coordinates  $r=\theta^0, \theta^1,\ldots, \theta^n$ on a neighborhood
   $(A,\Omega)\times \mathcal{U}$ about any point.
Observe that for each fixed $r$, equation (\ref{eq:rsolij}) 
is a tensorial identity on $S^n$, and that we may
therefore differentiate it using the Levi-Civita connection $\delt$ of $\gt$ to obtain
\begin{equation*}
	\delt_k\delt_i\delt_j f + \om\op f_{0k}\gt_{ij} = 0.
\end{equation*}
Hence
\[
	\delt_k\delt_i\delt_j f - \delt_i\delt_k\delt_j f
		= \om\op \left(f_{0i}\gt_{jk} - f_{0k}\gt_{ij}\right).
\]
On the other hand, since $\widetilde{R}_{ijkl} = \gt_{il}\gt_{jk}- \gt_{ik}\gt_{jl}$,
the standard commutation identities imply
\begin{eqnarray*}
	\delt_k\delt_i\delt_j f - \delt_i\delt_k\delt_j f &=& -\widetilde{R}_{kijl}\gt^{lm}f_m\\
		&=& -\left(\gt_{ij}\gt_{kl}-\gt_{il}\gt_{jk}\right)\gt^{lm}f_m\\
		&=& f_i \gt_{jk} -f_k\gt_{ij}  .
\end{eqnarray*}
Combining the two, we find
\begin{equation*}
	\left(\om\op f_{0i} - f_i\right)\gt_{jk} = \left(\om\op f_{0k}-f_k\right)\gt_{ij},
\end{equation*}
and tracing yields
\begin{equation*}
	(n-1)\left(\om\op f_{0i} - f_i\right) =0
\end{equation*}  
for all $i = 1,\ldots, n$.  Together with (\ref{eq:rsol0i}), we conclude
\begin{equation}
\label{eq:mainflat}
	(n-1)|\delt f|_{\gt}^2\left[1-(\op)^2\right] = 0.
\end{equation}

Since we assume $n > 1$, if $(Xf)(r_0,\mathbf{\theta_0})\neq 0$ for some $X\in T_{(r_0,\mathbf{\theta_0})}\Mnp$ tangent
to the $S^n$ factor, we must have $|\op| \equiv 1$ and $\opp\equiv 0$ on an interval $(a,b)\subset (A,\Omega)$ containing $r_0$.  
But, by equation (\ref{eq:sect}), this means that $\nu_1=\nu_2 = 0$ on $(a,b)$.  
We claim that $\nu_1= \nu_2 =0$ on the entire interval $(A,\Omega)$.

Let
\[
\beta = \sup \left\{r < \Omega\,|\, (\op)^2 = 1 \mbox{ on } (a, r)\right\}.
\]
If $\beta < \Omega$, by equations (\ref{eq:rsolij}) and (\ref{eq:mainflat}), we must have that $(\delt\delt f)(\beta,\cdot) = 0$,
$\op(\beta)=\sigma\in\{\pm 1\}$, $\opp(\beta) = 0$, $\om(\beta) > 0$, and $f_0(\beta,\cdot) = \sigma\lambda\om(\beta)$.
Moreover, for some small $\epsilon$,  $f$ is a function only of $r$
on $[\beta,\beta+\epsilon)$ and on this interval, $f^{\prime}$ and $\om$ satisfy the system (\ref{eq:rsys}), with the above initial conditions.
But one may check that the functions 
\[
\bar{\om}(r) = \sigma(r-\beta) + \om(\beta)\,\,\mbox{ and } 
	\bar{f}^{\prime}(r) = \lambda((r-\beta) +\om(\beta))
\] 
also satisfy (\ref{eq:rsys})
and agree with $\om$ and $f$ at $r=\beta$.  Therefore, by uniqueness
\footnote{Writing $x=\op$, and $u = f^\prime$, we may recast (\ref{eq:rsys}) as a first-order system
\begin{equation*}
	\begin{cases}
		\op &= x := F(\om, x, u)\\
		x^{\prime} &= xu -\lambda\om + (n-1)\frac{1-x^2}{\om}:=G(\om, x, u)\\
		u^{\prime} &= n\frac{xu}{\om} 
		+ (n-1)\left(n\frac{1-x^2}{\om^2}-\lambda\right):=H(\om, x, u).
	\end{cases}
\end{equation*}
Since $F$, $G$, $H$ are $C^{\infty}$ on the region $\{\om \neq 0\}$, the asserted uniqueness
follows from standard ODE theory.
},
these solutions must coincide and it follows that that $(\op)^2 =(\bar{\om}^{\prime})^2 = 1$ on the interval $[\beta,\beta +\epsilon)$,
contradicting our choice of $\beta$.  So $g$ is flat on $(a,\Omega)\times S^{n}$.  Using 
a similar argument at the other endpoint $a$, can show that $g$ must be flat on the entire cylinder $(A,\Omega)\times S^n$. 
 But this means either $\op\equiv 1$ or $\op\equiv -1$, so $g$ cannot extend to a smooth metric
on the sphere $S^{n+1}$ or to a complete metric on the cylinder.  The only possibility
is $\op\equiv 1$ and $\Mnp\approx \R{n+1}$.
\end{proof}

\subsection{An equivalent first-order system and its linearization.}

In view of the  result of the last section, we now assume $f=f(r)$. 
We are interested in solutions $(\om(r), f(r))$ to the system (\ref{eq:rsys})  for which 
$\om$ is strictly positive.

As Ivey observes in \cite{RFV2}, (\ref{eq:rsys}) is
invariant under translations of $r$ and $f$. 
By the introduction of the variables
\[
    x = \op\,\, \mbox{ and }\,\,  y=n\op-\om\fp
\]
which share this invariance, and of an independent variable $t$ which satisfies $dt = 1/\om dr$, one obtains
the first-order system
\begin{equation}
\label{eq:tsys}
	\begin{cases}
		\ddt{\om} &= x\om\\
		\ddt{x} &= x^2 -xy + n-1 -\lambda\om^2\\
		\ddt{y} &= xy - nx^2 -\lambda \om^2.
	\end{cases}
\end{equation}

Any solution to (\ref{eq:rsys}) gives rise to a trajectory of (\ref{eq:tsys}) and
conversely,
from a trajectory $(\om(t), x(t), y(t))$ of (\ref{eq:tsys}), one may recover $r$, $\om(r)$, and $f(r)$
by a succession of quadratures (see \cite{I}).  Consequently, 
it suffices to analyze solutions
to the simpler system (\ref{eq:tsys}).  
We take as coordinates $(\omega, x, y)$ on
the phase space $\R{3}$ and restrict our attention to trajectories lying in the 
half space $\om > 0$.

For $n > 1$, system (\ref{eq:tsys}) has two equilibrium points: 
$P_0 := (0, 1, n)$ and $P_1 :=  (0, -1, -n)$.  
Denoting the right hand side of (\ref{eq:tsys}) by $\Phi$, one finds
\begin{equation*}
  d\Phi_{P_0} = -d\Phi_{P_1} 
  = \left(\begin{array}{ccc} 1 & 0 & 0\\    
      0 & 2-n & -1\\
      0 & -n & 1\end{array}\right)
\end{equation*}         
which has eigenvalues $2$, $1$, and $1-n$.  
Since we assume $n\geq 2$, both $P_0$ and $P_1$ are saddle points: $P_0$ ($P_1$)
lying at the intersection of a two-dimensional unstable (stable) manifold and a one-dimensional stable (unstable)
manifold.  In particular, there is a one-parameter family of trajectories in the half-space ${\om > 0}$
initially tangent to $(1,0,0)$, among which, in light of Lemma 
(\ref{lem:smooth}), lie the trajectories which give rise
to smooth solutions on $S^{n+1}$ and $\R{n+1}$ (see Examples (\ref{ex:spheresol}) and (\ref{ex:flatsol}) below).  Trajectories
which correspond to smooth solutions on $S^{n+1}$ must, in addition, tend to $P_1$ as $t\to\infty$, and hence lie in the
intersection of the global unstable and stable manifolds of $P_0$ and $P_1$, respectively.

\begin{remark}
\label{rem:refl}
  If $L:\R{3}\to\R{3}$ is the map $(\om, x, y)\mapsto (\om, -x,-y)$, then from any solution $\gamma(t) =(\om(t), x(t), y(t))$ of (\ref{eq:tsys})
  on $(S, T)$,
  one may obtain a new solution 
\[
  \bar{\gamma}(t) := L\left(\gamma(\tau(t))\right) 
= \left(\om(\tau(t)), -x(\tau(t)), -y(\tau(t))\right)
\]
on an appropriate interval $(\bar{S}, \bar{T})$
where $\tau$ is chosen to satisfy $\ddt{\tau} = -1$ and $\tau(\bar{S}) = T$, 
$\tau(\bar{T}) = S$.

By use of this device, one immediately
obtains that the global stable and 
unstable manifolds of $P_i$, $S_i$ and $U_i$, $i=1,2$ are related by
$L(S_0) = U_1$, $L(U_0) = S_1$, and moreover, if a set $V\subset\R{3}$
is preserved by the system (\ref{eq:tsys}) for increasing $t$ (i.e,
$\gamma(t_0)\in V$ implies $\gamma(t)\in V$ for $t > t_0$, as long
as the solution is defined), then $L(V)$ is preserved by the system for
decreasing $t$.
\end{remark}

\subsection{The standard examples.}
\label{ss:ex}
At this juncture, it is worthwhile to recall the standard solutions 
to (\ref{eq:rsys}) and locate the corresponding solution 
to (\ref{eq:tsys}) in $\om x y$-space.  The content of Theorem \ref{thm:coverclass}
is that this list essentially exhausts the possibilities for complete
solutions.

\begin{example}{(Round sphere)}
\label{ex:spheresol}
	The condition $f\equiv \mathrm{const}$ in (\ref{eq:rsys})
 leads to the constant curvature
        soliton structure 
\begin{equation*}
\label{eq:rsphere}
\begin{cases}
	\om(r) &=\sqrt{\frac{n}{\lambda}}\sin
			\left(\sqrt{\frac{\lambda}{n}}r\right)\\
	f(r)  &=\mathrm{const}
\end{cases}
\end{equation*}
on the sphere $S^{n+1}$.  The corresponding trajectory in $\om x y$-space
is the elliptical arc $\{nx^2+\lambda\om^2 = n\}$ 
lying in the plane $\{y=nx\}$ joining $P_0$ and $P_1$. 
\end{example}
\begin{example}{(Gaussian soliton)}
\label{ex:flatsol}
        The condition $\op\equiv 1$ in (\ref{eq:rsys}) leads to the flat
soliton structure
\begin{equation}
	\begin{cases}
		\om(r) &= r\\
		f(r) &= \frac{\lambda}{2} r^2 + \mathrm{linear} 
	\end{cases}
\end{equation}
	on $\R{n+1}$ and corresponds to the trajectory  $y= n-\om^2$ in the plane $\{x=1\}$ emanating
        from $P_0$.  Applying the device of Remark (\ref{rem:refl}), one obtains a similar trajectory in the plane
	$\{x= -1\}$, satisfying $y = -n +\lambda\om^2$ and tending to $P_1$ as $t\to\infty$. 
\end{example}
\begin{example}
\label{ex:cylsol}{(Standard cylinder)}
	The condition $\op\equiv 0$ in (\ref{eq:rsys}) leads to the structure
\begin{equation*}
	\begin{cases}
		\omega(r) &= \sqrt{\frac{n-1}{\lambda}}\\
		f(r) &= \frac{\lambda}{2}r^2 + \mathrm{linear}
	\end{cases}
\end{equation*}
on $\R{}\times S^{n}$ corresponding to the line $\{(\sqrt{(n-1)/\lambda}, 0, y)\}$ in $\om x y$-space.
\end{example}
\begin{remark}
\label{rem:xconst}
   The flat and cylindrical solutions described in the previous two examples are the only for which $\op\equiv \mathrm{const}$, and
the corresponding trajectories in $\om x y$-space describe 
the intersections of the planes 
$\{x=1\}$, $\{x=0\}$, and $\{x=-1\}$ with the set
$\{x^2 -xy + n-1 -\lambda\om^2 = 0\}$.  This leads to an observation
which we shall have repeated occasion to use in the sequel:  the only solutions $\gamma(t)$ of (\ref{eq:tsys}) for which
$\ddt{x}(t_0) = 0$ and $x(t_0) = 1$, $0$, or $-1$ at some $t_0$ are, by uniqueness,
those for which $x(t)\equiv 1$, $0$, or $-1$,
respectively.
\end{remark}

Finally we note (as does Ivey in \cite{RFV2}) that by taking $\lambda = 0$ 
in the $x$ and $y$-components of (\ref{eq:tsys}) one recovers
the analogous system for rotationally symmetric solutions
to the steady soliton equation. Thus the 
trajectories of (\ref{eq:tsys})  which lie in the plane $\{\om=0\}$ 
are naturally associated with steady soliton structures (although, of course,
the warping function $\om(r)$ of these structures no longer corresponds 
directly to the $\om$-coordinate).
In particular, the (one-dimensional) intersection of the unstable manifold of $P_0$ with the plane $\{\om=0\}$ contains two candidates for a smooth
steady soliton on $\R{n+1}$: one with negative sectional curvature near the origin, 
which turns out to be incomplete,
and one with positive curvature near the origin, which is the well-known 
\emph{Bryant soliton}-- a complete steady soliton on $\R{n+1}$ of 
strictly positive curvature.   

\section{Proof of Theorem \ref{thm:coverclass}}

In what follows, $\gamma(t)=(\om(t),x(t), y(t))$ will represent
a trajectory of (\ref{eq:tsys}) defined for $t$ in what we may take
to be a \emph{maximal} interval $(S,T)$ with 
$-\infty \leq S < T \leq \infty$. To reduce the clutter of our expressions,
we shall usually suppress the dependence
of the components of the trajectory on the parameter $t$.

\subsection{Some invariant sets. }

We begin our analysis of trajectories of 
(\ref{eq:tsys}) by observing that the second and third derivatives of the $x$-component
have the following convenient expressions.

\begin{lemma}
\label{lem:xders}
	The $x$-component of any trajectory
	 $\gamma(t)$ of (\ref{eq:tsys})
	satisfies

\begin{align}
	\ddtt{x} 
\label{eq:xder2}
		 &= (n-1)x(x^2-1) + (3x-y)\ddt{x}
\end{align}
and 
\begin{align}
	\ddttt{x} 
\label{eq:xder3}
		  &= 2x\left[(2n-1)x - y\right]\ddt{x} 
		+ 2\left(\ddt{x}\right)^2 + (3x-y)\ddtt{x}.
\end{align}
\end{lemma}
\begin{proof}
	We compute
\begin{align*}
\ddtt{x} 
		 &= (2x-y)\ddt{x} -x^2y + nx^3 -\lambda\om^2x\\
		 &= (n-1)x(x^2-1) + (3x-y)\ddt{x},
\end{align*}
and
\begin{align*}
	\ddttt{x} 
&= \left[(n-1)(3x^2-1) + \ddt{x} -\ddt{y}\right]\ddt{x}
		+ 2\left(\ddt{x}\right)^2 + (3x-y)\ddtt{x}\\
		  &= 2x\left[(2n-1)x - y\right]\ddt{x} 
		+ 2\left(\ddt{x}\right)^2 + (3x-y)\ddtt{x}.
\end{align*}
\end{proof}

With the above expressions we may easily 
establish the following qualitative results on the
behavior of trajectories of the system.

\begin{lemma}
\label{lem:xpres}
	The regions
\begin{equation*}
	\left\{x \geq 1 \mbox{ ,  }\ddt{x} \geq 0 \right\}
		\mbox{ , } \left\{ x \leq -1 \mbox{ ,  }\ddt{x} \leq 0\right\}
	\mbox{, and } \left\{y \leq 0\right\}
\end{equation*}
	are preserved under system (\ref{eq:tsys}) for increasing $t$,
	and
\begin{equation*}
	\left\{x \geq 1 \mbox{ ,  }\ddt{x} \leq 0 \right\}
		\mbox{ , } \left\{ x \leq -1 \mbox{ ,  }\ddt{x} \geq 0\right\}
	\mbox{, and } \left\{y \geq 0\right\},
\end{equation*}
	are preserved for decreasing $t$.
\end{lemma}
\begin{proof}
	By (\ref{eq:xder2}), if $\ddt{x}(t_0) = 0$ at some $t_0$ with
	$x(t_0) > 1$, then $\ddtt{x}(t_0) > 0$, and hence both
	$\ddt{x}$ and $x$ continue to increase.  By Remark \ref{rem:xconst},
	we cannot have $\ddt{x} = 0$ and $x = 1$ simultaneously
	unless $x\equiv 1$.  The argument 
	for the case $x \leq -1$ and $\ddt{x} \leq 0$ is similar.  
	
        To see that $\{y\leq 0\}$ is preserved, observe that
\[
	\ddt{y} = -nx^2-\lambda\om^2 < 0
\]
	whenever $y = 0$.

	The preservation of the remaining sets for decreasing $t$
	follows by applying the results already obtained to the trajectory 
	$\bar{\gamma}(t)= L(\gamma(\tau))$ constructed as in Remark (\ref{rem:refl}).
\end{proof}

\begin{lemma}
\label{lem:xsign}$\phantom{b}$
\begin{enumerate}
\item  
\label{it:xs1}	If there exists $t_0\in (S,T)$ at which
	$x(t_0) = 0$, $y(t_0)\leq 0$, and $\ddt{x}(t_0) > 0$,
	then $x(t)$ (and $\ddt{x}$) increase until
	$\gamma(t)$ enters the region $\{x > 1\}$.
\item
\label{it:xs2b}
	If there exists $t_0\in (S,T)$ at which
	$x(t_0) = 0$, $y(t_0)\geq 0$, and $\ddt{x}(t_0) > 0$,
	then $x(t) < 0$ and  $\ddt{x}(t) > 0$ for all $t < t_0$ and there exists a $t_1 \leq t_0$
	such that $x(t) < -1$ for all $t <t_1$.
\end{enumerate}
In particular, in view of Lemma \ref{lem:xpres},
 if the trajectory $\gamma$ enters the region $\{x < 0\}$, either
 it remains there or eventually lies in the region $\{x > 1\}$.	
\end{lemma}
\begin{proof}
	In case (\ref{it:xs1}), we have $y(t) < 0$ for all $t> t_0$
	by Lemma \ref{lem:xpres},	
 	and $\ddtt{x}(t_0) \geq 0$ by equation (\ref{eq:xder2}).  
	Since by (\ref{eq:xder3}),
 $\ddttt{x} > 0$ on the region
\[
	\left\{ x \geq 0,\,\,\ddt{x} > 0,\,\, \ddtt{x} \geq 0\right\},
\]
we have $\ddtt{x}(t) > 0$ for all $t > t_0$.  Consequently,
$x(t) > \ddt{x}(t_0)(t-t_0)$.  As bounds on $x$ imply bounds on the derivatives
of $\om$ and $y$, the solution cannot expire while $0 <x< 1$. 
Since the interval $(S,T)$ is maximal, $x > 1$ eventually.

Case (\ref{it:xs2b}) then follows from case (\ref{it:xs1}) by 
considering again the trajectory $\bar{\gamma}(t)$
constructed according to the device in Remark \ref{rem:refl}:    
If $x(t_0)= 0$, $y(t_0) \geq 0$ and $\ddt{x}(t_0) > 0$, then $\bar{x}(t_0)=0$,
$\bar{y}(t_0)\leq 0$ and $\ddt{\bar{x}}(t_0) > 0$, and the corresponding interval
of definition $(\bar{S},\bar{T})$, satisfying $\tau(\bar{T}) = S$, $\tau(\bar{S}) = T$
will also be maximal. 

\end{proof}

\subsection{Proof of the case $\Mnp\approx S^{n+1}$}

\begin{theorem}
\label{thm:pco}
	Suppose $(S^{n+1}, g)$ 
	is a rotationally symmetric shrinking soliton.
	Then $g$ has positive curvature operator.
\end{theorem}
\begin{proof}
	Recall that, for a rotationally symmetric metric, the positivity of the
curvature operator is implied by that of the sectional curvatures, and these
	curvatures have the expressions
\[
	\nu_1 = -\frac{1}{\om^2}\ddt{x}\,\mbox{, and }\, \nu_2 = \frac{(1-x^2)}{\om^2}
\]
	in terms of the $(\om,x, y)$ coordinates. 
	
	Any trajectory $\gamma(t)$ of system (\ref{eq:tsys}) which corresponds to a smooth soliton structure on the sphere
	must tend to $P_0 = (0,1,n)$ as $t\to-\infty$ and to $P_1 = (0,-1, -n)$ as $t\to\infty$.
	Thus, by Lemma \ref{lem:xpres}, we must have $-1 < x(t) < 1$ (hence $\nu_2(t) > 0)$ for all $t$ 
and $\ddt{x} < 0$ at least initially. 
	We wish to show $\ddt{x} < 0$ for all $t$.

	 By equation (\ref{eq:xder2}) of Lemma \ref{lem:xders}, 
	 $\ddtt{x} < 0$ at critical points of $x$ in the region $\{0 < x < 1\}$,
	 so $\ddt{x}$ remains strictly negative on this region
         and, by Remark \ref{rem:xconst},
	 cannot vanish on $\{x=0\}$. 	
	Since $\gamma(t)$ must tend to 
	$P_1$ as $t\to\infty$, it must, in particular, enter
	the region $\{x < 0\}$, and, in view of Lemma \ref{lem:xsign}, 
	remain there for all subsequent $t$.
	Since $\ddtt{x} > 0$ at critical points of $x$ in the region 
	$\{-1 < x < 0\}$, $\ddt{x}$ must
	therefore remain strictly negative if $x$ is to approach $-1$. 
	So $\ddt{x} < 0$ always, and
	thus for any trajectory emanating from $P_0$ and tending to $P_1$ 
	we have $\nu_1 > 0$ and $\nu_2 > 0$ for all $t$ as claimed.
	
	Taking limits, one finds that at the ``poles'' $r = A$ and 
	$r  = \Omega$,
	the sectional curvatures agree and are at least non-negative.
	One may therefore apply Lemma 8.2 of \cite{H4} to conclude
	that the curvature operator $Rm(g):\wedge^2\to\wedge^2$ 
	is of constant rank, and therefore strictly positive everywhere.
\end{proof}

That $g$ has constant sectional curvature
is then an immediate consequence of the following

\begin{theorem*}{(B\"{o}hm and Wilking \cite{BW}: Theorem 1)}
On a compact manifold the normalized Ricci flow evolves a Riemannian
metric with $2$-positive curvature operator to a limit metric with
constant sectional curvature.
\end{theorem*}

That $f\equiv\mathrm{const}$ in 
case (\ref{it:sn}) of Proposition \ref{thm:coverclass}, 
then follows by substituting
$\nu_1 = \nu_2 = \mathrm{const}$ into (\ref{eq:rsys}) or, alternatively,
by considering the identity
\[
	R + \left|\nabla f\right|^2 -2\lambda f = \mathrm{const}
\]
valid on any gradient Ricci shrinking soliton (see, e.g., \cite{RFV2}).  
If  $f$ attains its maximum and minimum at the points $x_M$ and $x_m$, 
respectively, then the identity implies $f(x_M) =f(x_m)$ 
since $R$ is
constant.

\subsection{The asymptotic behavior of trajectories corresponding
	to complete, non-compact metrics}
\label{ss:asymp}
Hereafter, we shall consider solutions 
$\gamma(t)$ satisfying one or both of the conditions
 
\begin{equation}
\label{eq:compcond}
	\int_{t_0}^{T}\om(\sigma)\,d\sigma = \infty = \lim_{t\to T} r(t)
\end{equation}
and
\begin{equation}
\label{eq:compcond2}
	\int_{S}^{t_0}\om(\sigma)\,d\sigma = \infty = -\lim_{t\to S} r(t)
\end{equation}
for any $t_0\in (S,T)$.  Condition (\ref{eq:compcond})
is necessarily satisfied by any trajectory corresponding to a complete metric
on $\R{n+1}$ and both (\ref{eq:compcond}) and (\ref{eq:compcond2}) are necessarily
satisfied by any
trajectory corresponding to a complete metric on
 $\R{}\times S^n$.
As we shall see, these conditions impose rather stringent
conditions on the asymptotic behavior of a trajectory.

We remark that if $\gamma(t)$ satisfies (\ref{eq:compcond}) then 
$\bar{\gamma}(t) = L(\gamma(\tau))$ satisfies (\ref{eq:compcond2}) and vice-versa.  Thus,
from the following results, which apply to trajectories satisfying the ``forwards''
condition (\ref{eq:compcond}), we may easily obtain corresponding results for trajectories
satisfying the ``backwards'' condition (\ref{eq:compcond2}).  These results will be collected
in Lemma \ref{lem:banalog}, below.

We begin with a simple consequence of 
the forwards extendability condition by which we may obtain eventual
knowledge of the sign of the $y$-component.
\begin{lemma}
\label{lem:yneg}
	Along any trajectory $\gamma(t)$, the quantity $Q = y/\om$ 
	is strictly decreasing.  If $\gamma(t)$ satisfies (\ref{eq:compcond}), 
	$\lim_{t\to T} Q = -\infty$. In particular, $y$ eventually becomes negative.
\end{lemma}
\begin{proof}
	We compute
\[
	\dt Q = -n\frac{x^2}{\om} -\lambda\om < 0.
\]
	Integrating, we find that, for any $S <t_0 < t <T$,
\[
	\frac{y}{\om}(t) \leq \frac{y}{\om}(t_0)  -\lambda\int_{t_0}^t
		\om(\sigma) d\sigma.
\]
\end{proof}
The following observation is also immediate.
\begin{lemma}
\label{lem:xsup}
If $\gamma(t)$ satisfies (\ref{eq:compcond}), then
	$\limsup_{t\to T} x(t) \geq 0$. 
\end{lemma}
\begin{proof}
	If $x(t) < -\delta$  on some $(a, T)\subset (S, T)$, then
	$\om(t) \leq Ce^{-\delta t}$ on the same interval, and $\gamma(t)$ cannot
	satisfy (\ref{eq:compcond}). 
\end{proof}

Thus, in view of Lemma \ref{lem:xsign}, no trajectory satisfying (\ref{eq:compcond}) can enter the 
region $\{x < -1\}$. That no trajectory satisfying (\ref{eq:compcond}) can enter the region
$\{x > 1\}$ is true (as we prove next), but less obvious since $\om(t)\to\infty$.  
Taken together, Lemmas \ref{lem:xsign}, \ref{lem:xsup}, and \ref{lem:toofast} prove
that any complete metric on $\R{n+1}$ satisfies $\nu_2 \geq 0$.

	The proof follows the lines of an argument due to Bryant and
Ivey (c.f. \cite{RFV2})
 demonstrating the incompleteness of a similar 
trajectory of the steady soliton system.  In fact, if one regards
the trajectories of \eqref{eq:tsys} in the plane $\{\om= 0\}$ as 
trajectories of the the analogous system for steady solitons (cf. the
 remarks at the end of Section \ref{ss:ex}), then Ivey's argument
pertains to the trajectory in the plane emerging from $P_0$ 
in the direction opposite the Bryant soliton. The following
lemma then may be viewed as an extension of his finding
to the neighboring family of trajectories
in the unstable manifold $U_0$ which populate the sector
between the flat trajectory with $x\equiv 0$
and the plane $\{\om=0\}$.  These trajectories correspond to metrics 
of strictly negative
sectional curvature and are all incomplete. 
However, Ivey's argument does not
carry over directly, as, in the expression for $\ddt{x}$ in the shrinking
case, one has 
to contend with an additional term ($-\lambda\om^2$) of uncooperative sign.

\begin{lemma}
\label{lem:toofast}
	Suppose $x(t_0) > 1$ and $\ddt{x}(t_0) > 0$
	at some $t_0\in (S,T)$.  
	Then $\int_{t_0}^{T} \om(\sigma)\,d\sigma <\infty$.
\end{lemma}
\begin{proof}
We proceed by contradiction.
Suppose $\gamma(t)$ satisfies
(\ref{eq:compcond}). Then, by (\ref{eq:xder2}), $x(t) > 1$ and
 $\ddt{x}(t) > 0$ for $t > t_0$, and,  
by Lemma \ref{lem:yneg},
there is a $t_1\in (t_0, T)$ such that $y(t) < 0$ for all $t \geq t_1$.

Hence, by (\ref{eq:xder2}), we have
\begin{equation}
\label{eq:ddxconv}
	\ddtt{x} \geq (n-1)x(x^2 -1) + 3x\ddt{x} >  \frac{3}{2}\ddt{(x^2)}.
\end{equation}

Now, since the interval $(S, T)$ is assumed maximal, and since
bounds on $x$ imply bounds on the derivatives of $y$ and $\om$, if $T <\infty$
we must have $\limsup_{t\to T} |x(t)| = \lim_{t\to T} x = \infty$. On the other hand,
$x(t)$ is uniformly convex by \eqref{eq:ddxconv}, so even if $T = \infty$ 
we still have $\lim_{t\to T}x(t) = \infty$.
Returning to \eqref{eq:ddxconv} with this fact in hand, we find
\[
	\ddt{x} \geq \frac{5}{4} x^2 + 1
\] 
for all $t$ greater than some $t_2\geq t_1$. 
(The coefficient $5/4$ is chosen for convenience and could be replaced by $3/2 -\epsilon$ 
for any $\epsilon > 0$ -- below, we
merely require it to be greater than one.)
From this equation it follows that
$T < \infty$ 
and
\[
	\dt \arctan\left(\frac{\sqrt{5}}{2}x\right) 
	\geq \frac{\sqrt{5}}{2},
\]
which implies
\[
	x(t) \leq \frac{4}{5(T-t)}
\]
for $t$ sufficiently close to $T$.
Since $\ddt{}\log\om = x $, integrating and applying the above bound yields
\[
	\om(t) \leq \frac{C_2}{(T-t)^{\frac{4}{5}}} 
\]
for some constant $C_2$, contradicting (\ref{eq:compcond}).
\end{proof}

Together, Lemmas \ref{lem:xpres} and \ref{lem:toofast} allow us to restrict
our attention to trajectories which remain in the region 
$\{ -1 < x < 1\}$, and, consequently, to those with infinite
existence time $t\in (S, \infty)$ since a trajectory with $x$ bounded
cannot satisfy condition (\ref{eq:compcond}) on a interval bounded above.
Along such trajectories, $y$ becomes negative and 
Lemma \ref{lem:xsign} implies that eventually 
$x$ acquires a constant sign. As a consequence, we obtain 
the following refinement of Lemma \ref{lem:yneg}:
\begin{lemma}
\label{lem:yinf}
	If $\gamma(t)$  is a trajectory of (\ref{eq:tsys}) defined
	on $(S, \infty)$ satisfying condition (\ref{eq:compcond}) and
	$-1<x(t)<1$, then
	$\lim_{t\to\infty} y(t) = -\infty$.	
\end{lemma}
\begin{proof}
	Suppose $\limsup_{t\to\infty} y(t) \geq -M$ for some 
	$M > 1$, and choose $t_k\in (S,\infty)$ such that
	$t_k\nearrow\infty$ 
	and  $y(t_k) \geq -M$.  Then, by Lemma \ref{lem:yneg}, 
	$\lim_{k\to\infty} \om(t_k) = 0$. Since $x$ eventually acquires
	a constant sign, and $\dt{\om} = x\om$, we must have 
	$x(t) \leq 0$ eventually, and $\om(t)$ must tend to $0$ outright
	as $t\to\infty$.

	Now, equation (\ref{eq:xder2}) shows that $x(t)$ cannot
	attain a local maximum on $\{-1<x\leq 0\}$ unless
	$x = 0$, and we know $x = 0$ and $\ddt{x} = 0$ simultaneously
	only if $x\equiv 0$, in which case $dy = -(n-1)dt$.
	Otherwise, $x$ is eventually monotonic in $t$ and either
	 increases or decreases to a limit $\bar{x}\in [-1,0]$.
	By the remarks preceding this lemma, we cannot
	have $\bar{x} < 0$  if the trajectory is to 
	satisfy condition (\ref{eq:compcond}). So assume $\bar{x} = 0$,
	which implies $\ddt{x} > 0$ eventually.  
	Then, for any $\epsilon$, we can choose $t_{\epsilon}$ such that
	$t > t_{\epsilon}$ implies both $y(t) < 0$ and 
	$nx^2 + \lambda\om^2 < \epsilon$.  For such $t$, we have
\[
	\ddt{y} = xy - nx^2 -\lambda\om^2 > -\epsilon. 
\]

	Fixing $\epsilon < \frac{n-1}{2}$, we find, each for $k$,
\[
	y(t) > -M -\epsilon(t-t_k)
\]
and
\[
	\ddt{x}(t) = x^2 -xy + n-1 -\lambda\om^2 > \frac{n-1}{2}
		+x(t_k)(M + \left(\epsilon(t-t_k)\right),
\]
(where, in obtaining the last inequality,  we used that $x$ is monotonically
increasing).
 
Hence,
\[
	x(t)-x(t_k) > 
	\left(\frac{n-1}{2} + Mx(t_k)\right)(t-t_k) 
	+\epsilon\frac{x(t_k)}{2}(t-t_k)^2.
\]
For $k >> 0$, $Mx(t_k)  > -(n-1)/4$, so that the above (with the monotonicity
of $x(t)$) implies that there exists
$\delta = \delta(M,\epsilon, n) > 0$ and a subsequence $t_{k_j}\to \infty$
such that $x(t_{k_{j+1}}) > x(t_{k_{j}}) + \delta$ for all $j$.  
This contradicts
that $x\nearrow 0$, and proves $\limsup_{t\to\infty} y(t) = -\infty$.
\end{proof}

\subsection{Proof of the case $\Mnp\approx \R{n+1}$}

The results of the last section are enough to assemble
the \begin{proof}[Proof of Claim (\ref{it:rn}) of Theorem 
\ref{thm:coverclass}]

Since the underlying manifold $\Mnp$ is diffeomorphic
to $\R{n+1}$, the smooth extension of the metric to the origin 
$r=0$ requires $S= -\infty$ and our solution 
$\gamma(t)= (\om(t), y(t), \om(t))$ of (\ref{eq:tsys})
to satisfy $\lim_{t\to -\infty} \gamma(t) = P_0 = (0, 1, n)$. 

We claim first that if our trajectory is to satisfy condition 
(\ref{eq:compcond}), then $x\leq 1$.  For if ever $x > 1$, 
since $\lim_{t\to -\infty} x(t) = 1$,
we would have to have $\ddt{x}(t_0) > 0$ and $x(t_0) > 1$
at some earlier $t_0$. But then, by Lemmas
\ref{lem:xpres} and \ref{lem:toofast}, the $x$-component would blow-up too 
fast for $\gamma(t)$ to satisfy (\ref{eq:compcond}).  
So we must have $x \leq 1$ for all $t$.

Then, if ever $x=1$, we must also have $\ddt{x} =0$ at the same time, which,
as pointed out in Remark \ref{rem:xconst}, happens only if $x\equiv 1$--i.e.,
only if $\gamma(t)$ corresponds to the flat solution
of Example \ref{ex:flatsol}.  We claim that this
is the only trajectory emanating from $P_0$ which satisfies 
\eqref{eq:compcond}.

We may now assume that $x(t) < 1$ on our trajectory and that  
for some $t_0$ (hence all $t < t_0$), $\ddt{x}(t_0) < 0$ and 
$0 < x(t_0) < 1$. Since, by equation (\ref{eq:xder2}), $\ddtt{x}$ is
strictly negative at all critical points of $x$ in the region
$0 < x < 1$, there are two possibilities for our trajectory: either
\begin{enumerate}
	\item \label{it:poslim} 
	$x$ decreases monotonically to a limit $\bar{x}\in [0,1)$ as 
		$t\to\infty$, or
	\item \label{it:neg}
		$\gamma(t)$ enters the region $\{x\leq 0\}$ at some time 
		$t= t_1$.
\end{enumerate}
Knowing that $y\to -\infty$ 
 as $t\to\infty$ (in fact, just knowing that eventually
$y< 0$ suffices), 
we can dispose of case 
(\ref{it:poslim}) by observing that while $x \in [0,1]$,
\begin{equation}
\label{eq:eventneg}
	\ddtt{x} \leq (n-1)x(x^2-1) +\epsilon \ddt{x} < 0
\end{equation}
once $y < -\epsilon < 0$.  Hence $x$ eventually becomes negative.

Now, we also know from Lemma \ref{lem:xsup}
that $x$ cannot tend to a negative limit or become
strictly less than $-1$ if condition (\ref{eq:compcond}) is to be satisfied.  
Since $\ddtt{x}$ is strictly positive at critical points of $x$ in
the region $\{-1 < x < 0\}$, and since the only 
	trajectories with critical points of $x$ on the boundary of this
	region are classified in Examples \ref{ex:flatsol} and \ref{ex:cylsol}
and neither emanate from $P_0$, we face only two alternatives: 

\begin{enumerate}
\item[(2a)] either $\gamma$ enters the region $x > 0$ again, or
\item[(2b)]  $x\nearrow 0$ as $t\to\infty$.
\end{enumerate}  

Alternative (2a) is immediately excluded by Lemmas \ref{lem:xsign} and 
\ref{lem:toofast}: no trajectory which emanates from $P_0$ can
satisfy $x(t_0) = 0$, $\ddt{x}(t_0) > 0$, $y(t_0) \geq 0$, and
no trajectory which satisfies  $x(t_0) = 0$, $\ddt{x}(t_0) > 0$, 
$y(t_0) \leq 0$ can satisfy \eqref{eq:compcond}.

For (2b), we observe that since
$y\to -\infty$ as $t\to\infty$, we have $3x - y > \epsilon > 0$
eventually, and thus we may obtain the analog of equation (\ref{eq:eventneg})
for $t$ sufficiently large
\begin{equation}
\label{eq:eventpos}
	\ddtt{x} = (n-1)x(x^2-1) +(3x-y)\ddt{x} > \epsilon\ddt{x} >  0,
\end{equation}
which is incompatible with $x\nearrow 0$ as $t\to\infty$. 

The trajectory $x\equiv 1$ is therefore the unique trajectory
emanating from $P_0$ satisfying (\ref{eq:compcond}), and the
proof of the case $\Mnp\approx \R{n+1}$ is complete.
\end{proof}

\subsection{Proof of the case $\Mnp\approx \R{}\times S^{n}$}

As remarked earlier, the results in Section
\ref{ss:asymp}, regarding trajectories satisfying the
forwards extendability condition \eqref{eq:compcond}
have natural analogs for trajectories satisfying
the backwards version \eqref{eq:compcond}.

\begin{lemma}$\phantom{b}$
\label{lem:banalog}
\begin{enumerate}
\item
\label{it:e2}
	Suppose $\gamma(t)$ satisfies \eqref{eq:compcond2}.
\begin{enumerate}
\item $\lim_{t\to S} Q(t) = \infty$ and $y$ is initially positive.

\item 
\label{it:xinf} 
$\liminf_{t\to S} x(t) \leq 0$
\item If $-1 < x(t) < 1$ for all $t$ (so $-S=T=\infty$), 
then $\lim_{t\to -\infty}y(t)=\infty$.
\end{enumerate}
\item 
\label{it:btoofast}
If $x(t_0) < -1$ and $\ddt{x}(t_0) > 0$ at some $t_0\in (S,T)$,
\begin{equation}
\label{eq:bfin}
	\int_S^{t_0} \om(\sigma)\, d\sigma <\infty
\end{equation}
\end{enumerate}
\end{lemma}
\begin{proof}
	Let $\bar{\gamma}(t) = L(\gamma(\tau(t)))$ be
	 as in Remark \ref{rem:refl}.
\begin{equation}
\label{eq:rint}
	\int_{S}^{t_0}\om(\sigma)\,d\sigma 
= \int_{t_0}^{\bar{T}}\bar{\om}(\sigma)\,
	d\sigma.
\end{equation}
Thus, if $\gamma(t)$ satisfies \eqref{eq:compcond2}, $\bar{\gamma}(t)$
satisfies \eqref{eq:compcond} and the claims of part \ref{it:e2}
follow by the application of
 Lemmas \ref{lem:yneg}, \ref{lem:xsup}, and \ref{lem:yinf} to 
$\bar{\gamma}(t)$.

For part \ref{it:btoofast}, note that $\bar{x}(t_0) > 1$, $\ddt{\bar{x}}(t_0) > 0$
if $x(t_0) < -1$, $\ddt{x}(t_0) > 0$, so Lemma \ref{lem:toofast}
and equation \eqref{eq:rint} yield the inequality \eqref{eq:bfin}.
\end{proof}

Now we turn to the remainder of the proof of Theorem \ref{thm:coverclass}.

\begin{proof}[Proof of Claim \ref{it:cyl} of Theorem \ref{thm:coverclass}]

We shall show that the only trajectory $\gamma(t)$ satisfying 
both (\ref{eq:compcond}) and (\ref{eq:compcond2}) is that of 
Example \ref{ex:cylsol} with $x\equiv 0$ and 
$\om\equiv \sqrt{(n-1)/\lambda}$.

By Lemmas \ref{lem:xpres} and \ref{lem:toofast}
and case (\ref{it:btoofast}) of Lemma \ref{lem:banalog},
we can assume that $-1\leq x(t)\leq 1$ and 
consequently also $S = -\infty$, $T=\infty$.  
In fact, since 
neither of the trajectories
with $x\equiv \pm 1$ can satisfy both conditions (\ref{eq:compcond}) and (\ref{eq:compcond2}),
we may assume $ -1 < x(t) < 1$.

By Lemma \ref{lem:xsup} and part (\ref{it:xinf})
of Lemma \ref{lem:banalog}, we know 
$\liminf_{t\to -\infty} x(t)\leq 0$ and $\limsup_{t\to\infty} x(t) \geq 0$.
Thus since $\ddtt{x}$ is strictly positive at critical points of $x$ in the
region $\{-1 < x < 0\}$ and strictly negative at these critical points
in the region $\{0 < x < 1\}$, either
 $\gamma(t)$ crosses the plane $\{x= 0\}$ at some time
or its $x$-component maintains a constant sign and
tends monotonically to $0$ as $t\to\infty$ or $t\to -\infty$. 
We claim that this latter option cannot occur.  
For,  in light of the remarks in the preceding paragraph, the only     
scenarios of this option in which $\gamma(t)$ could potentially satisfy
the extendability criteria would be $x < 0$ and $x\nearrow 0$ as $t\to \infty$,
or $x > 0$ and $x\searrow 0$ as $t\to -\infty$.  
But the case $x < 0$, $x\nearrow 0$ was eliminated in the argument
for Claim (\ref{it:rn}), and the case $x\searrow 0$ as $t\to -\infty$
reduces to the previous one by the consideration of the trajectory 
$\bar{\gamma}(t)$.

Thus, we
conclude that there must exist a $t_0$ such that $x(t_0) = 0$.
If $\ddt{x}(t_0)\neq 0$, then we may assume $\ddt{x}(t_0) > 0$,
as the argument given in the case $\Mnp\approx \R{n+1}$ implies
that the only trajectories satisfying (\ref{eq:compcond}) and 
$\ddt{x}(t_0) < 0$ initially lie in the region $\{x < -1\}$ and cannot
thus satisfy (\ref{eq:compcond2}).
However, if $\ddt{x}(t_0) > 0$, then Lemma \ref{lem:xsign} implies that either
again $\gamma(t)$ lies initially in the region $\{x < -1\}$ or eventually
in the region $\{x > 1\}$, in which case, by Lemma \ref{lem:toofast}
and part (\ref{it:btoofast}) of Lemma \ref{lem:banalog}, $\gamma(t)$ can satisfy
at most one of the conditions (\ref{eq:compcond}) and (\ref{eq:compcond2}).

Thus we can only have
$\ddt{x}(t_0) =0$, which implies that $\gamma$ coincides with the 
trajectory $x\equiv 0$, $\om \equiv\sqrt{(n-1)/\lambda}$ as claimed.
\end{proof}
\begin{acknowledgement*}
	The author wishes to thank Professors Bennett Chow and Lei Ni
	for their support, encouragement, and 
	for useful discussions pertaining to this work. The author is 
	especially
	indebted to Professor Thomas Ivey, whose analysis of the problem
	in the expanding and steady soliton cases served as a guide for the study
	of the shrinking soliton case presented here.
\end{acknowledgement*}

\end{document}